\documentclass[a4paper]{article}
\usepackage{mathptmx}
\usepackage{amscd,amssymb,amsmath,amsthm}
\usepackage{mathrsfs}

\makeatletter
\renewcommand{\@seccntformat}[1]
{{\csname the#1\endcsname}.\hspace{0.3em}}

\renewcommand{\section}{\@startsection
{section}
{1}
{0mm}
{-1.5\baselineskip}
{\baselineskip}
{\bfseries\normalsize}}

\renewcommand{\subsection}{\@startsection
{subsection}
{2}
{0mm}
{-\baselineskip}
{0.5\baselineskip}
{\normalsize\itshape}}

\makeatother

\theoremstyle{plain}
\newtheorem{theorem}{Theorem}
\newtheorem{lemma}{Lemma}

\newtheorem{prop}{Proposition}

\newtheorem*{width1}{Width Inequality~I}
\newtheorem*{width2}{Width Inequality~II}
\newtheorem*{ML1}{Main Lemma~I}
\newtheorem*{ML2}{Main Lemma~II}

\theoremstyle{definition}
\newtheorem*{defin*}{Definition}

\theoremstyle{remark}

\newtheorem*{example*}{Example}
\newtheorem*{remark*}{Remark}

\DeclareMathAlphabet{\matheur}{U}{eur}{m}{n}
\DeclareMathAlphabet{\matheus}{U}{eus}{m}{n}
\DeclareMathAlphabet{\matheuf}{U}{euf}{m}{n}

\numberwithin{equation}{section}
\newcommand{\abs}[1]{\left\lvert#1\right\rvert}
\newcommand{\norm}[1]{\left\lVert#1\right\rVert}

\DeclareMathOperator{\iso}{Isom}

\DeclareMathOperator{\rank}{rank}
\begin{document}

\author{Gerasim  Kokarev %\thanks{Supported by EPSRC ???}
\\ {\small\it School of Mathematics, The University of Edinburgh}
\\ {\small\it King's Buildings, Mayfield Road, Edinburgh EH9 3JZ, UK}
\\ {\small\it Email: {\tt G.Kokarev@ed.ac.uk}}
}

\title{On geodesic homotopies of controlled width and conjugacies in
  isometry groups}
\date{}
\maketitle

\begin{abstract}
\noindent
We give an analytical proof of the Poincar\'e-type inequalities for
widths of geodesic homotopies between equivariant maps valued in
Hadamard metric spaces. As an application we obtain a linear bound for
the length of an element conjugating two finite lists in a group
acting on an Hadamard space.
\end{abstract}


%\noindent
%{\bf Keywords}: width of homotopy, harmonic maps, Hadamard space,
%decision problems.}

%\tableofcontents

\section{Introduction}
Let $M$ and $M'$ be smooth Riemannian manifolds without boundary. For
a smooth mapping $u:M\to M'$ by $E(u)$ we denote its energy
\begin{equation}
\label{e}
E(u)=\int\limits_M\norm{du(x)}^2d\mathit{Vol}(x),
\end{equation}
where the norm of the linear operator $du(x):T_xM\to T_{u(x)}M'$ is
induced by the Riemannian metrics on $M$ and $M'$. Let $u$ and $v$ be
smooth homotopic mappings of $M$ to $M'$ and $H(s,\cdot)$ be a smooth
homotopy between them. The {\it $L_2$-width} $W_2(H)$ of $H$ is
defined as the $L_2$-norm of the function
\begin{equation}
\label{length}
\ell_H(x)=\text{ the length of the curve }s\mapsto H(s,x),\qquad x\in M.
\end{equation}
A smooth homotopy $H(s,x)$ is called {\it geodesic} if for each $x\in
M$ the track curve $s\mapsto H(s,x)$ is a geodesic. 

In~\cite{KKS03,KKS05} Kappeler, Kuksin, and Schroeder prove the
following geometric inequality for the $L_2$-widths of geodesic
homotopies when the target manifold $M'$ is non-positively curved.
\begin{width1}
Let $M$ and $M'$ be compact Riemannian manifolds and suppose that $M'$
has non-positive sectional curvature. Let $\zeta$ be a homotopy class
of maps of $M$ to $M'$. Then there exist constants $C_\star$ and $C$
with the following property: any smooth homotopic maps $u$ and
$v\in\zeta$ can be joined by a geodesic homotopy $H$ whose $L_2$-width
$W_2(H)$ is controlled by the energies of $u$ and $v$,
\begin{equation}
\label{W1}
W_2(H)\leqslant C_\star(E^{1/2}(u)+E^{1/2}(v))+C.
\end{equation}
Moreover, if the sectional curvature of $M'$ is strictly negative, the
constants $C_\star$ and $C$ can be chosen to be independent of the
homotopy class $\zeta$.
\end{width1}

This  inequality can be viewed as a version of the Poincar\'e
inequality for mappings between manifolds. It also has an
isoperimetric flavour; it says that the `measure' of the cylinder
induced by the homotopy is estimated in terms of the `measure' of its
boundary. Inequality~\eqref{W1} is a key ingredient in the proof of
compactness results for perturbed harmonic map
equation~\cite{KKS03,GK1}. The latter, combined with old results of
Uhlenbeck, yields Morse inequalities for harmonic maps with
potential~\cite{GK2}.

The proof of Width Inequality~I in~\cite{KKS03,KKS05} is based on an
analogous inequality for maps of metric graphs;
see~\cite[Th.~5.1]{KKS03}. In more detail, let $G$ be a finite graph
and $u:G\to M'$ be a smooth map, that is whose restriction to every
edge is smooth. The length $L(u)$ of $u$ is defined as the sum of the
lengths of the images of the edges. By the $L_\infty$-width
$W_\infty(H)$ of a homotopy $H$ we mean the $L_\infty$-norm of the
length function $\ell_H(x)$, given by~\eqref{length}.
\begin{width2}
Let $G$ be a finite graph and $M'$ be a compact manifold of
non-positive sectional curvature. Let $\zeta$ be a homotopy class of
maps $G\to M'$. Then there exist constants $C_\star$ and $C$ with the
following property: any smooth homotopic maps $u$ and $v\in\zeta$ can
be joined by a geodesic homotopy $H$ such that
$$
W_\infty(H)\leqslant C_\star(L(u)+L(v))+C.
$$
Moreover, if the sectional curvature of $M'$ is strictly negative, the
constants $C_\star$ and $C$ can be chosen to be independent of the
homotopy class $\zeta$.
\end{width2}

The purpose of this note is two-fold: firstly, we generalise the width
inequalities to the framework of equivariant maps valued in Hadamard
spaces. This, in particular, includes width inequalities for
homotopies between maps into non-compact metric target spaces. In
contrast with the geometric methods in~\cite{KKS03,KKS05} (and
also in~\cite{BH}), we give an analytical proof of the width
inequalities via harmonic map theory. 

Secondly, we use width inequalities for equivariant maps of trees
to obtain informaion on algebraic properties of finitely generated
groups $\Lambda$ acting by isometries on Hadamard spaces. More
precisely, under some extra hypotheses, these groups satisfy the
following property: given two finite conjugate lists of elements
$(a_i)_{1\leqslant i\leqslant N}$ and $(b_i)_{1\leqslant i\leqslant
  N}$ in $\Lambda$ there exists $g\in\Lambda$ with $b_i=g^{-1}a_ig$
such that
$$
\abs{g}\leqslant C_\star\sum_{i=1}^N(\abs{a_i}+\abs{b_i})+C,
$$
where $\abs{\cdot}$ stands for the length $d(\cdot,e)$ in the word
metric on $\Lambda$. If the group $\Lambda$ has a soluble word
problem, then the latter estimate yields immediately the solubility of
the conjugacy problem for finite lists in $\Lambda$.

\section{Statements and discussion of results}
\subsection{Width inequalities for equivariant maps}
Let $M$ be a compact Riemannian manifold without boundary; we denote
by $\tilde M$ its universal cover and by $\Gamma$ the fundamental
group $\pi_1(M)$. Let $(Y,d)$ be an Hadamard space; that is a complete
length space of non-positive curvature in the sense of Alexandrov (see
Sect.~\ref{prem} for a precise definition). Denote by $\rho$ a
representation of $\Gamma$ in the isometry group of $Y$. Recall that a
map $u:\tilde M\to Y$ is called $\rho$-{\it equivariant} if
$$
u(g\cdot x)=\rho(g)\cdot u(x)\quad\text{for all}\quad x\in\tilde M,
\ g\in\Gamma. 
$$
For $\rho$-equivariant maps $u$ and $v$ the real-valued functions
$d(u(x),v(x))$, where $x\in\tilde M$, are invariant with respect to
the domain action and, hence, are defined on the quotient $M=\tilde
M/\Gamma$. In particular, the quantity
\begin{equation}
\label{L2d}
d_2(u,v)=\left(\int_M d^2(u(x),v(x))\mathit{dVol}(x)\right)^{1/2}
\end{equation}
defines a metric on the space of locally $L_2$-integrable
$\rho$-equivariant maps. The latter can be also regarded as the
$L_2$-width of a unique geodesic homotopy between $\rho$-equivariant
maps. If $u$ is a locally Sobolev $W^{1,2}$-smooth $\rho$-equivariant
map, then its energy density measure $\abs{du}^2\mathit{dVol}$ (see
Sect.~\ref{prem}) is also $\Gamma$-invariant and the energy of $u$ is
defined as the integral
\begin{equation}
\label{E}
E(u)=\int_M\abs{du}^2\mathit{dVol}.
\end{equation}

Recall that the {\it ideal boundary} of $Y$ is defined as the set of
equivalence classes of {\it asymptotic} geodesic rays, where two rays
are asymptotic if they remain at a bounded distance from each
other. Clearly, any action of $\Gamma$ by isometries on $Y$ extends to
the action on the ideal boundary.

\begin{theorem}
\label{t1}
Let $M$ be a compact Riemannian manifold without boundary and $Y$ be a
locally compact Hadamard space. Let $\Gamma$ be the fundamental group
of $M$ and $\rho:\Gamma\to\iso(Y)$ be its representation whose image does
not fix a point on the ideal boundary of $Y$. Then there exists a
constant $C_\star$ such that for any $\rho$-equivariant locally
$W^{1,2}$-smooth maps $u$ and $v$ the $L_2$-width of a geodesic
homotopy $H$ between them satisfies the inequality
\begin{equation}
\label{L2width}
W_2(H)\leqslant C_\star(E^{1/2}(u)+E^{1/2}(v)).
\end{equation}
\end{theorem}
The proof appears in Sect.~\ref{proofs1}. The idea is to prove first a
similar inequality when one of the maps is an energy minimiser, and
then to use compactness properties of the moduli space formed by such
minimisers. The former is based on a compactness argument,
mimicking the proof of the classical Poincar\'e inequality.

Below we state a version of Theorem~\ref{t1} for equivariant maps of
trees. First, we introduce more notation.
Let $G$ be a finite connected graph without terminals and $\Gamma$ be
its fundamental group $\pi_1(G)$. By $T$ we denote the universal
covering tree of $G$; the group $\Gamma$ acts naturally on $T$ by the
deck transformations. As above the symbol $\rho$ denotes a
representation of $\Gamma$ in the isometry group of an Hadamard space
$Y$. For a locally rectifiable $\rho$-equivariant map $u:T\to Y$, its
length density measure $\abs{du}dt$ (see Sect.~\ref{prem}) is
$\Gamma$-invariant and the length of $u$ is defined as the integral
$$
L(u)=\int_G\abs{du}dt.
$$
\begin{theorem}
\label{t1g}
Let $G$ be a finite graph and $Y$ be a locally compact Hadamard
space. Let $\Gamma$ be the fundamental group of $G$ and $\rho:
\Gamma\to\iso (Y)$ be its representation whose image does not fix a
point on the ideal boundary of $Y$. Then there exists a constant
$C_\star$ such that for any locally rectifiable $\rho$-equivariant
maps $u$ and $v$ the $L_\infty$-width of a geodesic homotopy between
them satisfies the inequality
\begin{equation}
\label{L_infty_width}
W_\infty(H)\leqslant C_\star(L(u)+L(v)).
\end{equation}
\end{theorem}
\begin{example*}
Let $M'$ be a (not necessarily compact) Riemannian manifold whose
sectional curvature is negative and bounded away from zero and the
injectivity radius is positive. Let $\rho:\Gamma\to\pi_1(M')$ be a
homomorphism whose image is neither trivial nor infinite cyclic. Then
the latter does not fix a point on the ideal boundary of the universal
cover of $M'$. Indeed, the group $\rho(\Gamma)$ is generated by
hyperbolic elements (regarded as isometries of the universal cover),
see~\cite[Lem.~B.1]{KKS05}, and the statement follows from the results
in~\cite[Sect.~6]{EO}. Thus, as a particular case, Theorem~\ref{t1g}
contains the width inequality for homotopies between maps from $G$ to
$M'$; this is the situtation considered
in~\cite[Th.~0.1]{KKS05}. (Under the hypotheses on the homomorphism
$\rho$, the homotopy class is neither trivial nor contains a map onto
a closed curve.) The methods in~\cite{KKS03,KKS05} do not seem to
yield an analogous $L_2$-width inequality (provided by
Theorem~\ref{t1}) for non-compact targets when the dimension of the
domain is greater than one.
\end{example*}

We proceed with width inequalities for representations in
co-compact subgroups of $\iso(Y)$. Recall that an action of a group
$\Lambda$ on a metric space $(Y,d)$ is said to be {\it co-compact} if
the quotient $Y/\Lambda$ is compact. Further, the action of $\Lambda$
is said to be {\it proper} if for each $y\in Y$ there exists $r>0$
such that the set $\{g\in\Lambda~|~g\cdot B(y,r)\cap
B(y,r)\ne\varnothing \}$ is finite. For a homomorphism
$\rho:\Gamma\to\Lambda$ we denote by $Z$ below the centraliser of the
image $\rho(\Gamma)$ in $\Lambda$.
\begin{theorem}
\label{t2}
Let $M$ be a compact Riemannian manifold without boundary and $Y$ be a
locally compact Hadamard space. Let $\Lambda$ be a group acting
properly and co-compactly by isometries on $Y$. Denote by $\Gamma$ the
fundamental group of $M$ and let $\rho$ be a homomorphism $\Gamma\to
\Lambda$. Then there are constants $C_\star$ and $C$ such that for any
$\rho$-equivariant locally $W^{1,2}$-smooth maps $u$ and $v$ there
exists an element $h\in Z$ such that the $L_2$-width of a geodesic
homotopy $H$ between $u$ and $h\cdot v$ satisfies the inequality
$$
W_2(H)\leqslant C_\star(E^{1/2}(u)+E^{1/2}(v))+C.
$$
\end{theorem}
\begin{theorem}
\label{t2g}
Let $G$ be a finite graph and $Y$ be a locally compact Hadamard
space. Let $\Lambda$ be a group acting properly and co-compactly by
isometries on $Y$. Denote by $\Gamma$ the fundamental group of $G$ and
let $\rho$ be a homomorphism $\Gamma\to\Lambda$. Then there are
constants $C_\star$ and $C$ such that for any locally rectifiable
$\rho$-equivariant maps $u$ and $v$ there exists an element $h\in Z$
such that the $L_\infty$-width of a geodesic homotopy $H$ between
$u$ and $h\cdot v$ satisfies the inequality
$$
W_\infty(H)\leqslant C_\star(L(u)+L(v))+C.
$$
\end{theorem}
\begin{remark*}
If the homomorphism $\rho:\Gamma\to\Lambda$ in the theorems is
trivial, then the second constant $C$ is equal to $\mathit{diam}
(Y/\Lambda)\mathit{Vol}^{1/2}M$ and $\mathit{diam}(Y/\Lambda)$ in
the $L_2$- and $L_\infty$-versions respectively. For non-trivial
representations of $\Gamma$ it can be chosen to be zero. 
\end{remark*}
\begin{example*}
As a partial case, when the action of $\Lambda$ is free,
Theorems~\ref{t2} and~\ref{t2g} above contain width inequalities for
homotopies between continuous $W^{1,2}$-smooth maps valued in a
compact metric space $Y/\Lambda$. The choice of an element $h\in Z$ in
this setting corresponds to the choice of the homotopy between
maps. Indeed, recall that the fundamental group of the space formed by
continuous maps homotopic to  $u:M\to Y/\Lambda$ is equal to the
centraliser of the image $u_*(\pi_1(M))$ in $\Lambda$.
\end{example*}

\subsection{Conjugacies of finite lists in isometry groups}
Now we describe some applications of the width inequalities to
geometric group theory. First, recall that a discrete subgroup
$\Lambda$ in a Lie group $\mathbf G$ is called {\it lattice} if the
quotient $\mathbf G/\Lambda$ carries a finite $\mathbf G$-invariant
measure. Such a lattice is always finitely generated provided the
group $\mathbf G$ is semi-simple and has $\rank\geqslant 2$; see
ref. in~\cite{LMR}. Choose a finite system of generators $(g_i)$ of
$\Lambda$ and consider the word metric $d(\cdot,\cdot)$ on $\Lambda$
associated with the Cayley graph determined by the generators. Denote
by $\abs{g}$ the length $d(g,e)$, the distance between an element $g$
and the neutral element $e$.

The following statements are essentially consequences of
Theorems~\ref{t1g} and~\ref{t2g} and are explained in
Sect.~\ref{groups}.

\begin{theorem}
\label{group1}
Let $\mathbf G$ be a semi-simple Lie group of $\rank\geqslant 2$ all
of whose simple factors are non-compact. Let $\Lambda$ be an
irreducible lattice in $\mathbf G$ and $(a_i)_{1\leqslant i\leqslant
  N}$ be a finite list of elements in $\Lambda$ which does not fix a
point on the ideal boundary of the associated symmetric space. Then
for any conjugate (in $\Lambda$) list $(b_i)_{1\leqslant i\leqslant
  N}$ any conjugating element $g\in\Lambda$, $b_i=g^{-1}a_ig$,
satisfies the inequality
$$
\abs{g}\leqslant C_\star\sum_{i=1}^N(\abs{a_i}+\abs{b_i})+C,
$$
where the constants depend only on the conjugacy class of the lists.
In particular, for such two given lists the set of conjugating
elements is finite. 
\end{theorem}
\begin{remark*}
An analogous statement holds if $\Lambda$ is an irreducible lattice in
an almost simple $p$-adic algebraic Lie group of $\rank\geqslant
2$. In this case we consider lists which do not fix points on the
ideal boundary of the associated Euclidean building.
\end{remark*}
\begin{example*}
When the group $\mathbf G$ is algebraic, the hypothesis on the finite
list $(a_i)$ is satisfied if, for example, the elements $a_i$'s
generate a lattice (e.g., the whole group $\Lambda$) in $\mathbf
G$. Indeed, by Borel's density theorem the latter is Zariski dense in
$\mathbf G$ and, hence, does not fix a point on the ideal boundary of
the associated symmetric space.
\end{example*}

The estimate above yields immediately an algorithm deciding whether a
given list of elements in $\Lambda$ is conjugate to the list $(a_i)$
in the theorem. This is a special case of the more general result due
to Grunewald and Segal~\cite{GS}: the conjugacy problem for finite
lists in arithmetic groups is soluble. (Any irreducible lattice in a
semi-simple Lie group of $\rank\geqslant 2$ is arithmetic, by the
Margulis theorem.) However, we do not know whether the linear estimate
for the length of the conjugating element holds under weaker
hypotheses than in Theorem~\ref{group1}.

We proceed with the conjugacy problem for finite lists in groups which
act properly and co-compactly on Hadamard spaces by isometries. Recall
that such groups are necessarily finitely presented;
see~\cite[I.8.11]{BH99}. As above by $\abs{g}$ we denote the length
$d(g,e)$ in the word metric.
\begin{theorem}
\label{group2}
Let $Y$ be a locally compact Hadamard space and $\Lambda$ be a group
acting properly and co-compactly by isometries on $Y$. Then for any
finite conjugate lists $(a_i)_{1\leqslant i\leqslant N}$ and 
$(b_i)_{1\leqslant i\leqslant N}$ of elements in $\Lambda$ there
exists an element $g\in\Lambda$ with $b_i=g^{-1}a_ig$ such that
\begin{equation}
\label{group_ineq}
\abs{g}\leqslant C_\star\sum_{i=1}^N(\abs{a_i}+\abs{b_i})+C,
\end{equation}
where the constants depend only on the conjugacy class of the lists.
Further, there exists an algorithm deciding whether two given finite
lists of elements in $\Lambda$ are conjugate.
\end{theorem}

When the list $(a_i)$ in the theorem consists of a single element, the
solubility of the conjugacy problem is well-known. It is, for example,
a consequence of an exponential (compare with our linear) bound for
the length of the conjugating element in~\cite[III.$\Gamma$.1.12]{BH99}.
In the context of decision problems it is worth noting that there
are finitely presented groups in which the conjugacy problem for
elements is soluble, but the conjugacy problem for finite lists is
not. We refer to~\cite{BH} for the explicit examples. Finally, mention
that in~\cite{BH} Bridson and Howie prove a closely related linear
estimate for the length of the conjugating (two finite lists) element
in Gromov hyperbolic groups.

\section{Preliminaries}
\label{prem}
\subsection{Sobolev spaces of maps to metric targets}
We recall some background material on Sobolev spaces of maps
valued in a metric space. The details can be found in~\cite{KS}.

Let $\Omega$ be a Riemannian domain and $(Y,d)$ be an arbitrary metric
space. We suppose that $\Omega$ is endowed with a Lebesgue measure
$\mathit{dVol}$ induced by the Riemannian volume. A measurable map
$u:\Omega\to Y$ is called {\it locally $L_2$-integrable} if it has a
seperable essential range and for which $d(u(\cdot),Q)$ is a locally
$L_2$-integrable function on $\Omega$ for some $Q\in Y$ (and, hence,
by the triangle inequality for any $Q\in Y$). If the domain $\Omega$
is bounded, then the function 
$$
d_2(u,v)=\left(\int_\Omega d^2(u(x),v(x))\mathit{dVol}(x)
\right)^{1/2}
$$
defines a metric on the space of locally $L_2$-integrable maps. The
latter is complete provided $Y$ is complete.

The {\it approximate energy density} of a locally $L_2$-integrable map
$u$ is defined for $\varepsilon>0$ as
$$
e_\varepsilon(u)(x)=\int_{S_\varepsilon(x)}\frac{d^2(u(x),u(x'))}
{\varepsilon^{n+1}}\mathit{dVol}(x'),
$$
where $S_\varepsilon(x)$ denotes the $\varepsilon$-sphere centred at
$x$ and $n$ stands for the dimension of $\Omega$. The function
$e_\varepsilon(x)$ is non-negative and locally $L_1$-integrable.
\begin{defin*}
The {\it energy} $E(u)$ of a locally $L_2$-integrable map $u$ is
defined as
$$
E(u)=\sup_{0\leqslant f\leqslant 1}\left(\lim_{\varepsilon\to 0}
\sup\int_\Omega fe_\varepsilon(u)\mathit{dVol}\right),
$$
where the sup is taken with respect to compactly supported continuous
functions which take values between $0$ and $1$. A locally
$L_2$-integrable map $u$ is called {\it locally $W^{1,2}$-smooth} if
for any relatively compact domain $D\subset\Omega$ the energy
$E(\left.u\right|_D)$ is finite.
\end{defin*}
Due to the results of Korevaar and Schoen~\cite[Sect.~1]{KS} a locally
$L_2$-integrable map $u$ is locally $W^{1,2}$-smooth if and only if 
there exists a locally $L_1$-integrable function $e(u)$ such that the
measures $e_\varepsilon(u)\mathit{dVol}$ converge weakly to the
measure $e(u)\mathit{dVol}$ as $\varepsilon\to 0$. The function
$e(u)$, also denoted by $\abs{du}^2$, is called the {\it energy
  density} of $u$, and the energy $E(u)$ is equal to the total mass
$\int e(u)\mathit{dVol}$.

Now suppose that the domain $\Omega$ is $1$-dimensional, that is an
interval $I=(a,b)$. For a map $u:I\to Y$ one can also define the {\it
  approximate length density} as
$$
l_\varepsilon(u)(t)=\frac{d(u(t),u(t+\varepsilon))+d(u(t),
u(t-\varepsilon))}{\varepsilon},\qquad t\in I.
$$
Then the length of $u$ is defined by the formula similar to that for
the energy,
$$
L(u)=\sup_{0\leqslant f\leqslant 1}\left(\lim_{\varepsilon\to 0}
\sup\int_I fl_\varepsilon(u)dt\right),
$$
where the sup is taken with respect to compactly supported continuous
functions. A map $u:I\to Y$ is called {\it rectifiable} if its length
is finite. In this case there exists a {\it length density function} 
(or {\it speed function}) $l(u)$ such that the lenght $L(u)$ equals
$\int l(u)dt$.

\subsection{Hadamard spaces}
Recall that an {\it Hadamard space} $(Y,d)$ is a complete metric space
which satisfies the following two hypotheses:
\begin{itemize}
\item[(i)] {\it Length Space.} For any two points $y_0$ and $y_1\in Y$
  there exists a rectifiable curve $\gamma$ from $y_0$ to $y_1$ such
  that 
$$
d(y_0,y_1)=\mathit{Length}(\gamma).
$$
We call such a curve $\gamma$ {\it geodesic}.
\item[(ii)] {\it Triangle comparison.} For any three points $P$, $Q$,
  and $R$ in $Y$ and the choices of geodesics $\gamma_{PQ}$,
  $\gamma_{QR}$, and $\gamma_{RP}$ connecting the respecting points
  denote by $\bar P$, $\bar Q$, and $\bar R$ the vertices of the
  (possibly degenerate) Euclidean triangle with side lengths
  $\ell(\gamma_{PQ})$, $\ell(\gamma_{QR})$, and $\ell(\gamma_{RP})$
  respectively. Let $Q_\lambda$ be a point on the geodesic
  $\gamma_{QR}$ which is a fraction $\lambda$,
  $0\leqslant\lambda\leqslant 1$, of the distance from $Q$ to $R$;
$$
d(Q_\lambda,Q)=\lambda d(Q,R),\qquad d(Q_\lambda,R)=(1-\lambda)d(Q,R).
$$
Denote by $\bar Q_\lambda$ an analogous point on the side $\bar Q\bar
R$ of the Euclidean triangle. The triangle comparison hypothesis says
that the metric distance $d(P,Q_\lambda)$ (from $Q_\lambda$ to the
opposite vertex) is bounded above by the Euclidean distance $\abs{\bar
P-\bar Q_\lambda}$. This inequality can be written in the following
form:
\begin{equation}
\label{H0}
d^2_{PQ_\lambda}\leqslant (1-\lambda)d^2_{PQ}+\lambda d^2_{PR}-
\lambda(1-\lambda)d^2_{QR}.
\end{equation}
\end{itemize}

It is a direct consequence of the property~(ii) above that geodesics
in an Hadamard space are unique. It is also a consequence of geodesic
uniqueness that an Hadamard space has to be
simply-connected~\cite[II.1]{BH99}. Examples include symmetric spaces
of non-compact type and Euclidean buildings, simply-connected manifolds
of non-positive sectional curvature, Hilbert spaces, simply-connected
Euclidean or hyperbolic simplicial complexes satisfying certain local
link conditions~\cite[II.5.4]{BH99}. Another class of examples is
provided by the following proposition.
\begin{prop}
\label{L2maps}
Let $M$ be a compact Riemannian manifold without boundary and $(Y,d)$
be an Hadamard space. Let $\rho$ be a represenation of the fundamental
group $\Gamma=\pi_1(M)$ in the group of isometries of $Y$. Then the
space of $\rho$-equivariant locally $L_2$-integrable maps from $\tilde
M$ to $Y$ endowed with the metric~\eqref{L2d} is an Hadamard space.
\end{prop}
The proof follows straightforward from the definitions: the geodesics
in the new space are geodesic homotopies and the triangle comparison
hypothesis follows by integration of relation~\eqref{H0}.

A useful consequence of the triangle comparison hypothesis is the
following quadrilateral comparison property due to
Reshetnyak~\cite{Re} (we refer to~\cite[Cor. 2.1.3]{KS} for a proof).
\begin{prop}
Let $(Y,d)$ be an Hadamard space and $P$, $Q$, $R$, and $S$ be an
ordered sequence of points in $Y$. For $0\leqslant\lambda,\mu\leqslant
1$ define $P_\lambda$ to be the point which is the fraction $\lambda$
of the way from $P$ to $S$ (on the geodesic $\gamma_{PS}$) and $Q_\mu$
to be the point which is the fraction $\mu$ of the way from $Q$ to $R$
(on the opposite geodesic $\gamma_{QR}$). Then for any
$0\leqslant\alpha,t\leqslant 1$ the following inequality holds:
\begin{equation}
\label{H1}
d^2_{P_tQ_t}\leqslant (1-t)d^2_{PQ}+td^2_{RS}-t(1-t)\left[\alpha
(d_{PS}-d_{QR})^2+(1-\alpha)(d_{RS}-d_{PQ})^2\right].
\end{equation}
%\begin{multline}
%\label{H2}
%d^2_{Q_tP}+d^2_{Q_{1-t}S}\leqslant d^2_{PQ}+d^2_{RS}+
%t(d^2_{SP}-d^2_{QR})+2t^2d^2_{QR}\\
%-t\left[\alpha(d_{SP}-d_{QR})^2+(1-\alpha)
%(d_{RS}-d_{PQ})^2\right].
%\end{multline}
\end{prop}
\noindent
Setting $\alpha$ to be equal to zero in this inequality, we
deduce the convexity of the distance between geodesics
\begin{equation}
\label{H3}
d_{P_tQ_t}\leqslant (1-t)d_{PQ}+td_{RS}.
\end{equation}
This implies the following energy convexity property. Let $u$ and $v$
be locally $W^{1,2}$-smooth maps from the Riemannian domain $\Omega$
to an Hadamard space $(Y,d)$. Let $H(s,\cdot)$ be a geodesic
homotopy between $u$ and $v$; the point $H(s,x)$ is the fraction $s$
of the way from $u(x)$ to $v(x)$, where $x\in\Omega$. Then for any $s$
the map $H(s,\cdot)$ is locally $W^{1,2}$-smooth and for any
relatively compact domain $D\subset\Omega$ its energy satisfies the
inequality
\begin{equation}
\label{Econvex}
E^{1/2}(H_s)\leqslant (1-s)E^{1/2}(u)+sE^{1/2}(v).
\end{equation}
Inequality~\eqref{H3} also yields the length convexity along geodesic
homotopies. More precisely, let $u$ and $v$ be rectifiable paths in
$(Y,d)$ and let $H(s,\cdot)$ be a geodesic homotopy between them
parameterised by the arc-length as above. Then for any $s$ the map
$H(s,\cdot)$ is rectifiable and its length satisfies the inequality
\begin{equation}
\label{Lconvex}
L(H_s)\leqslant (1-s)L(u)+sL(v).
\end{equation}

Another consequence of the triangle comparison hypothesis is the
existence of the {\it nearest point projection} $\pi:Y\to A$ onto a
convex subset $A$. In more detail, if $(Y,d)$ is an Hadamard space and
$A$ is its non-empty closed convex subset, then for any $y\in Y$ there
exists a unique point $a\in A$ which minimises the distance $d(y,a)$
among all points in $A$; see~\cite[Prop.~2.5.4]{KS}.

\subsection{Some properties of harmonic maps}
Let $M$ be a compact Riemannian manifold without boundary and $(Y,d)$
be an Hadamard space. As above by $\Gamma$ we denote the fundamental
group of $M$ and by $\rho:\Gamma\to\iso(Y)$ its representation in
the isometry group of $Y$. We consider $\rho$-equivariant locally
$W^{1,2}$-smooth maps $u$ from the universal cover $\tilde M$ to $Y$.
The energy density of such a map $u$ is a $\Gamma$-invariant function
on $\tilde M$, which can be also regarded as a function on the
quotient $M=\tilde M/\Gamma$. In particular, by the energy $E(u)$ we
understand the integral $\int_Me(u)\mathit{dVol}$.
We call a $\rho$-equivariant map {\it harmonic} if it minimises the
energy among all $\rho$-equivariant locally $W^{1,2}$-smooth maps.

The following statement is a straightforward consequence of the energy
convexity, formula~\eqref{Econvex}. We state it as a proposition for
the convenience of references.
\begin{prop}
\label{har1}
Under the hypotheses above, let $u$ and $v$ be two $\rho$-equivariant
harmonic maps and $H(s,\cdot)$ be a geodesic homotopy between them;
the point $H(s,x)$ is the fraction $s$ of the way from $u(x)$ and
$v(x)$, where $x\in\tilde M$. Then for each $s$ the map $H(s,\cdot)$
is also $\rho$-equivariant harmonic and the energy $E(H_s)$ does not
depend on $s$. 
\end{prop}

We proceed with the Lipschitz continuity of harmonic maps. The
following proposition is a consequence of the result by Korevaar and
Schoen~\cite[Th.~2.4.6]{KS}.
\begin{prop}
\label{har2}
Under the hypotheses above, any $\rho$-equivariant harmonic map $u$ is
Lipschitz continuous and its Lipschitz constant is bounded above by
$C\cdot E^{1/2}(u)$, where the constant $C$ depends on the manifold $M$
and its metric only.
\end{prop}

Now let $G$ be a finite connected graph without terminals and $\Gamma$
be its fundamental group. By $T$ we denote the universal covering tree
of $G$. Similarly to the discussion above, for a locally rectifiable
$\rho$-equivariant map $u:T\to Y$ the length density function $l(u)$
is $\Gamma$-invariant and, hence, descends to the quotient
$G=T/\Gamma$. In particular, by the length $L(u)$ we understand the
integral $\int_G l(u)dt$.
It is straightforward to see that if a map $u$ minimises the length
among all locally rectifiable $\rho$-equivariant maps, then its
restriction to every edge is a geodesic. If the latter has
constant-speed parameterisation on every edge, then it is also
harmonic and the length of every edge $u_I$ satisfies the relation
$L^2(u)=E(u_I)(b-a)$, see~\cite[Lemm.~12.5]{EeFu}. Conversely, if $u$
is a $\rho$-equivariant harmonic map, then its restriction to every
edge is a constant-speed geodesic whose squared length is proportional
to the energy as above. In particular, the length is constant on the
set of $\rho$-equivariant harmonic maps, where it achieves its
minimum.

\section{Proofs of the width inequalities}
\label{proofs1}
We start with the following lemma.
\begin{ML1}
Let $M$ be a compact Riemannian manifold without boundary and $(Y,d)$
be a locally compact Hadamard space. Let $\rho:\Gamma\to\iso(Y)$ be a
representation of the fundamental group $\Gamma=\pi_1(M)$. Suppose
that the moduli space $\mathit{Harm}$, formed by $\rho$-equivariant
harmonic maps, is non-empty and bounded in $L_2$-metric. Then there
exists a positive constant $C_\star$ with the following property: for
any $\rho$-equivariant locally $W^{1,2}$-smooth map $u$ there exists a
harmonic map $\bar u\in\mathit{Harm}$ such that
\begin{equation}
\label{MLineq}
d_2(u,\bar u)\leqslant C_\star(E^{1/2}(u)-E^{1/2}_\star),
\end{equation}
where $E_\star=E(\bar u)$ is the energy minimum among
$\rho$-equivariant maps.
\end{ML1}
\begin{proof}
First, note that inequality~\eqref{MLineq} is invariant under the
rescaling of the metric on the target space $Y$. Hence, it is
sufficient to prove the lemma under the assumption that
\begin{equation}
\label{scale}
\text{the distance }d_2(u,\bar u)\text{ is not less than one.}
\end{equation}
Suppose the contrary. Then there exists a sequence of maps $u_k$ such
that for any $\bar u\in\mathit{Harm}$
$$
d_2(u_k,\bar u)\geqslant k(E^{1/2}(u_k)-E^{1/2}(\bar u)).
$$
For each $u_k$ choose a harmonic map $\bar u_k$ at which the infimum
$$
d_2(u_k,\bar u_k)=\inf\{d_2(u,v): v\in\mathit{Harm}\}
$$
is attained. Such a harmonic map clearly exists: it is the value of
$u_k$ under the nearest point projection onto $\mathit{Harm}$.
(The lower semicontinuity of the energy~\cite[Th.~1.6.1]{KS} and
Prop.~\ref{har1} imply that $\mathit{Harm}$ is a closed convex subset
in the Hadamard space of $\rho$-equivariant locally $L_2$-integrable
maps.)

Denote by $H^k_s$, where $s\in [0,1]$, a geodesic homotopy between
$\bar u_k$ and $u_k$; we set $H^k_0=\bar u_k$ and
$H^k_1=u_k$. Assuming that the parameter $s$ is proportional to the
arc length, we obtain
$$
d_2(H^k_s,H^k_0)=s\cdot d_2(u_k,\bar u_k)\geqslant
s\cdot k(E^{1/2}(u_k)-E^{1/2}(\bar u_k)).
$$
Recall the energy $E^{1/2}(\cdot)$ is convex along geodesic
homotopies; 
$$
s(E^{1/2}(u_k)-E^{1/2}(\bar u_k))\geqslant
E^{1/2}(H^k_s)-E^{1/2}(H^k_0).
$$
Combining the last two inequalities we conclude that
\begin{equation}
\label{c}
d_2(H^k_s,H^k_0)\geqslant k(E^{1/2}(H^k_s)-E^{1/2}(H^k_0)).
\end{equation}
Now choose a sequence of $s_k\in [0,1]$ such that the distance
$d_2(H^k_{s_k},H^k_0)$ equals one; by the assumption~\eqref{scale}
this is possible. Then relation~\eqref{c} implies that the sequence
$E(H^k_{s_k})$ converges to $E_\star$ as $k\to +\infty$. Since the
moduli space $\mathit{Harm}$ is bounded in $L_2$-metric, the latter
together with the choice of the $s_k$'s implies that the sequence
$H^k_{s_k}$ is bounded in the $W^{1,2}$-sense; that is 
\begin{equation}
\label{star}
d_2(H^k_{s_k},w)+E(H^k_{s_k})\leqslant C,
\end{equation}
where $w$ is a fixed $\rho$-equivariant map. Now by the version of
Rellich's embedding theorem~\cite[Th.~1.13]{KS} we can find a
subsequence $H^k_{s_k}$ (denoted by the same symbol) which converges
in $L_2$-metric and point-wise to a locally $W^{1,2}$-smooth map $\bar
v$.
%$$
%d_2(H^k_{s_k},\bar v)\to 0\quad\text{as}\quad k\to +\infty.
%$$
By the lower semi-continuity of the energy~\cite[Th.~1.6.1]{KS} the map
$\bar v$ is energy minimising and by the point-wise convergence is
$\rho$-equivariant. By the choice of the $s_k$'s we clearly have
$$
d_2(H^k_1,H^k_{s_k})=d_2(H^k_1,H^k_0)-d_2(H^k_{s_k},H^k_0)=
d_2(u_k,\bar u_k)-1.
$$
Thus, the $L_2$-distance between the maps $u_k$ and $v$ can be
estimated as
$$
d_2(u_k,\bar v)\leqslant d_2(H^k_1,H^k_{s_k})+d_2(H^k_{s_k},\bar v)
=d_2(u_k,\bar u_k)+(d_2(H^k_{s_k},\bar v)-1).
$$
For sufficiently large $k$ the second term on the right-hand side is
negative and we arrive at a contradiction with the choice of the
harmonic maps $\bar u_k$'s.
\end{proof}

The following lemma summarises known results (essentially due
to~\cite{KS}) on the moduli space $\mathit{Harm}$, formed by
$\rho$-equivariant maps. 
\begin{lemma}
\label{harm1}
Let $M$ be a compact Riemannian manifold without boundary and $Y$ be a
locally compact Hadamard space. Let $\Gamma$ be the fundamental group
of $M$ and $\rho:\Gamma\to\iso(Y)$ be its representation whose image
does not fix a point on the ideal boundary of $Y$. Then the moduli
space $\mathit{Harm}$, formed by $\rho$-equivariant harmonic maps, is
non-empty and compact in $C^0$-topology.
\end{lemma}
Since there is no direct reference for the statement on the compactness
of $\mathit{Harm}$ and to make our paper more self-contained, we give
a proof now.
\begin{proof}[Proof of Lemma~\ref{harm1}]
First, we explain the existence of a $\rho$-equivariant harmonic
map. By~\cite[Th.~2.6.4]{KS} there exists an energy minimising
sequence $\{u_i\}$ of equivariant Lipschitz continuous maps, whose
Lipschitz constants are uniformly bounded. Let $\Omega$ be a
fundamental domain for the action of $\Gamma$ on the universal cover
$\tilde M$. We claim that under the hypotheses of the theorem the
ranges $u_i(\Omega)$ are contained in a bounded subset of $Y$. Indeed,
suppose the contrary. Then there exists a point $x\in\Omega$ such that
the sequence $\{u_i(x)\}$ is unbounded in $Y$, i.e.
$$
d(u_i(x),Q)\to +\infty\quad\text{for some}\quad Q\in Y.
$$
For any $g\in\Gamma$ consider the sequence $d(\rho(g)\cdot
u_i(x),u_i(x))$. By the equivariance of the $u_i$'s and the 
uniform boundedness of their Lipschitz constants we have
$$
d(\rho(g)\cdot u_i(x),u_i(x))\leqslant Cd(g\cdot x,x),
$$
and hence the quantities on the left hand side remain bounded as $i\to
+\infty$. By the convexity of the distance between geodesics,
relation~\eqref{H3}, we see that the (Hausdorff) distances between the
geodesic segments $\overline{Qu_i(x)}$ and $\rho(g)\cdot
\overline{Qu_i(x)}$ also remain bounded as $i\to +\infty$. Since $Y$
is locally compact, we can find a subsequence of $u_i$, denoted by the
same symbol, such that the segments $\overline{Qu_i(x)}$ converge on
compact subsets to a geodesic ray $\sigma$ with initial point at $Q$.
Then the distance between $\sigma$ and $\rho(g)\cdot\sigma$ is also
bounded for any $g\in\Gamma$. This shows that $\sigma$
represents a fixed point for the action of $\rho(\Gamma)$ and leads to
a contradiction.

Now, since $Y$ is locally compact, the Arzela-Ascoli theorem applies
and we can find a subsequence of $u_i$ converging in $C^0$-topology to
an energy-minimising and, hence, harmonic map. Thus, the moduli space
$\mathit{Harm}$ is non-empty. 

Finally, we explain the compactness of $\mathit{Harm}$. Let ${u_i}$ be
a sequence of $\rho$-equivariant harmonic maps. By Prop.~\ref{har1}
their energies coincide and Prop.~\ref{har2} the $u_i$'s are uniformly
Lipschitz continuous. The same argument as above shows that the ranges
$u_i(\Omega)$ are contained in a bounded subset of $Y$. Again by the
Arzela-Ascoli theorem there exists a converging subsequence. By the
lower semi-continuity of the energy the limit map is energy minimising
and, hence, harmonic. Thus, the moduli space $\mathit{Harm}$ is
compact in $C^0$-topology among $\rho$-equivariant harmonic maps.
\end{proof}
\begin{proof}[Proof of Theorem~\ref{t1}]
By Lemma~\ref{harm1} Main Lemma~I applies: for given
$\rho$-equivariant maps $u$ and $v$ we can find harmonic
$\rho$-equivariant maps $\bar u$ and $\bar v$ such that $d_2(u,\bar
u)$ and $d_2(v,\bar v)$ are estimated as in~\eqref{MLineq}. By
Lemma~\ref{harm1} the moduli space $\mathit{Harm}$ is compact and,
hence,  the distance between $d_2(\bar u,\bar v)$ is uniformly
bounded. The $L_2$-width of a geodesic homotopy $H$ between $u$ and
$v$ is the distance $d_2(u,v)$, and by the triangle inequality we have
$$
W_2(H)\leqslant d_2(u,\bar u)+d_2(\bar u,\bar v)+d_2(v,\bar v).
$$
The second term is bounded, and the first and the last can be
estimated as in~\eqref{MLineq}; thus, we obtain
$$
W_2(H)\leqslant C_\star(E^{1/2}(u)+E^{1/2}(v))+C,
$$
where $C$ equals $\mathit{diam}(\mathit{Harm})-2C_\star
E_\star^{1/2}$. Since, under the hypotheses of the theorem, the energy
minimum $E_\star$ is positive, this inequality can be re-written in
the form~\eqref{L2width}. 
\end{proof}

Now we explain the proof of Theorem~\ref{t1g}; it follows essentially
the same idea. First, we discuss the version of Main Lemma~I. By
$d_\infty(u,v)$ we denote below the maximum of the distance function
between maps $u$ and $v$.
\begin{ML2}
Let $G$ be a finite graph and $(Y,d)$ be a locally compact Hadamard
space. Let $\rho:\Gamma\to\iso(Y)$ be a representation of the
fundamental group $\Gamma=\pi_1(G)$. Suppose that the moduli space
$\mathit{Harm}$, formed by $\rho$-equivariant harmonic maps $T\to Y$,
is non-empty and compact in $C^0$-topology. Then there exists a
positive constant $C_\star$ with the following property: for any
continuous rectifiable $\rho$-equivariant  map $u$ there exists a
harmonic map $\bar u\in\mathit{Harm}$ such that
\begin{equation}
\label{MLineq2}
d_\infty(u,\bar u)\leqslant C_\star(L(u)-L_\star),
\end{equation}
where $L_\star=L(\bar u)$ is the length minimum among
$\rho$-equivariant maps.
\end{ML2}
\begin{proof}
First, without loss of generality we may assume that the maps $u:T\to
Y$ under consideration are such that their restrictions to every edge
are parameterised proportionally to the arc-length. Second, as in the
proof of Main Lemma~I, it is sufficient to prove the lemma under the
assumption that the distance $d_\infty(u,\bar u)$ is not less than
one.

Suppose the contrary. Then there exists a sequence of maps $u_k$ and
harmonic maps $\bar u_k$ such that
$$
d_\infty(u_k,{\bar u}_k)\geqslant k(L(u_k)-L_\star);
$$
we suppose that the ${\bar u}_k$'s minimise the distance
$\{d_\infty(u_k,\bar u)$, where $u\in\mathit{Harm}\}$. Denote by
$H^k_s$, where $s\in [0,1]$, a geodesic homotopy between $\bar u_k$
and $u_k$. Assuming that the parameter is proportional to the
arc-length and using the convexity of the length,
relation~\eqref{Lconvex}, we obtain
$$
d_\infty(H^k_s)\geqslant k(L(H^k_s)-L(H^k_0)).
$$
Choosing a sequence $s_k\in [0,1]$ such that the left-hand side above
equals to one, we conclude that $L(H^k_{s_k})$ converges to $L_\star$
as $k\to +\infty$. Since the lengths of $H^k_{s_k}$ are bounded and
the edges of the $H^k_{s_k}$'s are parameterised proportionally to the
arc-length, we see that the sequence of the $H^k_{s_k}$'s is
equicontinuous. Further, the compactness of $\mathit{Harm}$ implies
that the latter sequence is $d_\infty$-bounded. Now the Arzela-Ascoli
theorem applies and there exists a subsequence converging in
$d_\infty$-metric to a continuous map $\bar v$. The map $\bar v$ is
clearly $\rho$-equivariant and length-minimising. Moreover, it has a
constant-speed parametrisation and, hence, is harmonic. Now one gets
a contradiction in the same way as in the proof of Main Lemma~I.
\end{proof}
\begin{proof}[Proof of Theorem~\ref{t1g}]
First, Lemma~\ref{harm1} carries over the case of $\rho$-equivariant
maps of trees. In more detail, we need to start with a length
minimising sequence which is uniformly Lipschitz continuous. The
latter can be constructed by re-parameterising any length minimising
sequence proportionally to the arc-length on every edge. The rest of
the proof (of Lemma~\ref{harm1}) carries over without essential
changes.

Now we simply follow the lines in the proof of Theorem~\ref{t1} and
use Main Lemma~II instead of Main Lemma~I.
\end{proof}

We proceed with the proofs of Theorems~\ref{t2} and~\ref{t2g}. First,
recall some notation. Let $\Lambda$ be a group acting properly and
co-compactly by isometries on $Y$. For a homomorphism $\rho:\Gamma
\to\Lambda$ by $Z$ we denote the centraliser of the image
$\rho(\Gamma)$ in $\Lambda$. The group $Z$ acts naturally on the space
of $\rho$-equivariant maps $u:\tilde M\to Y$ and, in particular, on
the moduli space $\mathit{Harm}$.
\begin{lemma}
\label{harm2}
Under the hypotheses of Theorem~\ref{t2}, the moduli space
$\mathit{Harm}$, formed by $\rho$-equivariant harmonic maps, is 
non-empty and the quotient $\mathit{Harm}/Z$ is compact in
$C^0$-topology.
\end{lemma}
\begin{proof}
We start with the existence of a $\rho$-equivariant harmonic
map. By~\cite[Th.~2.6.4.]{KS} there exists an energy minimising
sequence $\{u_i\}$ of equivariant Lipschitz continuous maps, whose
Lipschitz constants are uniformly bounded. Let $\Omega$ and $D$ be 
 fundamental domains for the actions of $\Gamma$ on $\tilde M$ and
$\Lambda$ on $Y$ respectively. Fix a point $x_*\in\Omega$. Then there
exists a sequence of elements $h_i\in\Lambda$ such that the maps
$h_i\cdot u_i$ send $x_*$ into the closure of $D$. Since the $h_i$'s
are isometries, the sequence $\{h_i\cdot u_i\}$ is also energy
minimising and uniformly Lipschitz continuous. Moreover, since
$\Lambda$ acts co-compactly, its fundamental domain $D$ is bounded,
and the uniform Lipschitz continuity implies that the ranges $h_i\cdot
u_i(\Omega)$ are contained in a bounded subset of $Y$. By the
Arzela-Ascoli theorem there exists a subsequence, also denoted
by $h_i\cdot u_i$, converging to a limit map $v$.

Now we define a homomorphism $\varphi:\Gamma\to\Lambda$ such that the
limit map $v$ is $\varphi$-equivariant. For this fix a generator
$g\in\Gamma$ and consider the points
$$
v(g\cdot x)=\lim (h_i\cdot u_i)(g\cdot x)\quad\text{and}\quad
v(x)=\lim (h_i\cdot u_i)(x),
$$
where $x\in\Omega$. The triangle inequality implies that
$$
(h_i\rho(g)h_i^{-1})\cdot v(x)\to v(g\cdot x)\quad\text{as}
\quad i\to+\infty.
$$
Now, since the action of $\Lambda$ is proper, the sequence
$h_i\rho(g)h_i^{-1}$ contains a constant subsequence; we denote it
value by $\varphi(g)\in\Lambda$. We use the $h_i$'s of this subsequence
for the same procedure for another generator in $\Gamma$. Repeating
the process we define $\varphi$ on all generators. It then extends as
a homomorphism $\varphi:\Gamma\to\Lambda$ and the map $v$ is
$\varphi$-equivariant. As a result of this procedure, we also have
a sequence $h_i\in\Lambda$ such that
$$
h_i\rho(g)h_i^{-1}=\varphi(g)\quad\text{for any}\quad g\in\Gamma.
$$
This identity implies that the $h_i$'s can be written in the form
$k\cdot\bar{h_i}$, where $\bar{h_i}\in Z$, and the element
$k\in\Lambda$ conjugates $\rho$ and $\varphi$. Now, since the sequence
$h_i\cdot u_i$ converges to $v$, the sequence $\bar{h_i}\cdot u_i$
converges to $k^{-1}v$. Moreover, the latter is energy minimising and
is formed by $\rho$-equivariant maps. Thus, the limit map $k^{-1}v$ is
a harmonic $\rho$-equivariant map and the existence is demonstrated.

The compactness of $\mathit{Harm}/Z$ follows by the same argument as
above with the substitution of the sequence of harmonic maps for the
energy minimising sequence $\{u_i\}$. By Prop.~\ref{har1} the former
sequence is also energy minimising, and by Prop.~\ref{har2} is
uniformly Lipschitz continuous; the argument above yields a sequence
$\bar{h_i}\in Z$ such that $\bar{h_i}\cdot u_i$ converges to a
$\rho$-equivariant harmonic map.
\end{proof}

\begin{proof}[Proof of Theorem~\ref{t2}]
Let $\mathcal H$ be a fundamental domain for the action of $Z$ on the
moduli space $\mathit{Harm}$. First, Main Lemma~I holds under a weaker
hypothesis than the $L_2$-boundedness of $\mathit{Harm}$. More
precisely, it is sufficient to assume that the domain $\mathcal
H$ is bounded in the $L_2$-metric. Indeed, since the group $Z$ acts by
isometries, one can suppose that the maps $\bar u_k$'s (in the proof
of Main Lemma~I) belong to $\mathcal H$. The boundedness of the latter
is then used to obtain the $W^{1,2}$-boundedness of the sequence
$H^k_{s_k}$, relation~\eqref{star}. The rest of the proof stays
unchanged.

Now the combination of Lemma~\ref{harm2} and estimate~\eqref{MLineq}
yields the statement in the fashion similar to the proof of
Theorem~\ref{t1}.
\end{proof}

\begin{proof}[Proof of Theorem~\ref{t2g}]
First, Main Lemma~II holds under a weaker hypothesis than the
compactness of the moduli space $\mathit{Harm}$. Similarly to the
above, it is sufficient to assume that a fundamental domain for the
action of $Z$ on $\mathit{Harm}$ is compact. Further,
Lemma~\ref{harm2} carries over the case of $\rho$-equivariant maps of
trees; the proof follows essentially the same line of argument. The
combination of this version of Lemma~\ref{harm2} with
estimate~\eqref{MLineq2} yields the statement in the same fashion as
above.
\end{proof}

\section{Finitely generated subgroups in isometry groups}
\label{groups}

Recall that the action of a group $\Lambda$ on a metric space $(Y,d)$
by isometries defines an orbit pseudo-metric on $\Lambda$:
$$
d_y(g,h)=d(g\cdot y,h\cdot y),\quad\text{where}\quad g,h\in\Lambda,
$$
and $y\in Y$ is a fixed reference point. For another point $\bar y\in
Y$ the pseudo-metrics $d_y$ and $d_{\bar y}$ are coarsely isometric;
that is there exists a constant $C$ ($=2d(y,\bar y)$) such that
$$
d_{\bar y}(g,h)-C\leqslant d_y(g,h)\leqslant d_{\bar y}(g,h)+C.
$$
First, we show that the $L_\infty$-width inequalities imply an
estimate for the conjugating element in the orbit pseudo-metric.
\begin{lemma}
\label{orbit1}
Let $\mathbf{G}$ be a semi-simple Lie group all of whose simple
factors are non-compact. Let $\Lambda$ be an irreducible lattice in
$\mathbf{G}$ and $(a_i)_{1\leqslant i\leqslant N}$ be a finite list of
elements in $\Lambda$ which does not fix a point on the ideal boundary
of the associated symmetric space. Then for any conjugate (in
$\Lambda$) list $(b_i)_{1\leqslant i\leqslant N}$ any conjugating
element $g\in\Lambda$, $b_i=g^{-1}a_ig$, satisfies the inequality 
$$
d_y(g,e)\leqslant C_\star\sum_{i=1}^{N}(d_y(a_i,e)+d_y(b_i,e)), 
$$
where $y\in Y$ is a reference point, and the constant depends only on
the conjugacy class of the list $(a_i)$.
\end{lemma}
\begin{proof}
Let $Y$ be a symmetric space associated with the Lie group
$\mathbf{G}$. Under the hypotheses on $\mathbf{G}$, the natural
$\mathbf{G}$-invariant Riemannian metric on $Y$ defines a distance $d$
which makes $Y$ into an Hadamard space.

Consider the bouqet of $N$ copies of a circle; denote by
$\Gamma=\oplus_{i=1}^{N}\mathbf Z$ its fundamental group and by $T$
its universal cover. Define a homomorphism $\rho:\Gamma\to\Lambda$ by
the rule: the generator of the $i$th copy of $\mathbf{Z}$ maps into
$a_i$. For a fixed reference point $y\in Y$ consider the graph in $Y$
whose vertices are points $g\cdot y$, where $g$ is a word in the
alphabet $(a_i)$. The edges are geodesic arcs; two points $g_1\cdot y$
and $g_2\cdot y$ are joined by an edge if and only if $g_1^{-1}g_2$ is
an element $a_i$ or its inverse. Suppose that each edge is
parameterised proportionally to the arc-length. Such a parametrisation
defines a $\rho$-equivariant map $u:T\to Y$, whose length $L(u)$ is
given by the sum $\sum_{i=1}^Nd(a_iy,y)$.

Analogously, for a conjugate list $(b_i)_{1\leqslant i\leqslant N}$
one defines a $(g^{-1}\rho g)$-equivariant map $v:T\to Y$, where $g$ is
a conjugating element. Note that the map $g\cdot v$ is
$\rho$-equivariant and its length $L(g\cdot v)$ coincides with $L(v)=
\sum_{i=1}^Nd(b_iy,y)$. By the hypotheses of the lemma,
Theorem~\ref{t1g} applies and we have
$$
d(g\cdot y,y)\leqslant W_\infty(H)\leqslant C_\star(L(u)+L(v)),
$$
where $H$ is a homotopy between $u$ and $g\cdot v$. Now the
combination with the expressions for the lengths finishes the proof. 
\end{proof}
\begin{proof}[Proof of Theorem~\ref{group1}]
The statement is a direct consequence of Lemma~\ref{orbit1} and  the
solution of Kazhdan's conjecture in~\cite{LMR}. The latter says that
the word metric (with respect to some finite set of generators) on an
irreducible lattice $\Lambda$ is quasi-isometric to the orbit metric
(with respect to the action on the associated symmetric space or
Euclidean building) provided $\mathbf{G}$ is semi-simple and its
$\rank\geqslant 2$.
\end{proof}
\begin{lemma}
\label{orbit2}
Let $Y$ be a locally compact Hadamard space and $\Lambda$ be a group
acting properly and co-compactly by isometries on $Y$. Then for any
finite conjugate lists $(a_i)_{1\leqslant i\leqslant N}$ and 
$(b_i)_{1\leqslant i\leqslant N}$ of elements in $\Lambda$ there
exists an element $g\in\Lambda$ with $b_i=g^{-1}a_ig$ such that
$$
d_y(g,e)\leqslant C_\star\sum_{i=1}^{N}(d_y(a_i,e)+d_y(b_i,e))+C, 
$$
where $y\in Y$ is a reference point, and the constants depend only on
the conjugacy class of the list $(a_i)$.
\end{lemma}
\begin{proof}
The proof follows the same line of argument as the proof of
Lemma~\ref{orbit1} with the use of Theorem~\ref{t2g} instead of
Theorem~\ref{t1g}.
\end{proof}
\begin{proof}[Proof of Theorem~\ref{group2}]
By \v{S}varc-Milnor lemma~\cite[I.8.19]{BH99} the word and orbit metrics
on $\Lambda$ are quasi-isometric. The combination of this with
Lemma~\ref{orbit2} implies the first statement of the
theorem. Further, by~\cite[III.$\Gamma$.1.4]{BH99} the word problem in
$\Lambda$ is soluble. This yields the algorithm deciding the conjugacy
of finite lists in the following fashion. If there exists an element
conjugating two given lists, then it belongs to the finite subset of
$\Lambda$ formed by elements satisfying the
bound~\eqref{group_ineq}. Using the solubility of the word problem,
the algorithm checks all elements from this finite set.
\end{proof}

\end{document}